\documentstyle[multicol,prl,aps,epsf,psfig,epsfig]{revtex}
\def\be{\begin{equation}}
\def\ee{\end{equation}}
\begin{document}
\draft
\title{A Small World Network of Prime Numbers}
\author{Anjan Kumar Chandra and Subinay Dasgupta}

\address{Department of Physics, University of Calcutta,
92 A.P. C. Road, Calcutta 700009, India.\\
}
\maketitle
\begin{abstract}
According to Goldbach conjecture, any even number can be broken up as the
sum of two prime numbers : $n = p + q$. We construct a network where each
node is a prime number and corresponding to every even number $n$, we put
a link between the component primes $p$ and $q$. In most cases, an even
number can be broken up in many ways, and then we chose {\em one}
decomposition with a probability $|p - q|^{\alpha}$. Through computation
of average shortest distance and clustering coefficient, we conclude that for
$\alpha > -1.8$ the network is of small world type and for $\alpha < -1.8$
it is of regular type. 
We also present a theoretical justification for such behaviour.
\end{abstract}

\begin{multicols}{2}
\bigskip

\section{Introduction}

The study of networks, with an emphasis on small-world behaviour and scale 
invariant properties has turned out to be very important for analysing the
statistical properties of diverse type of systems \cite{review}. A network
is defined as a graph consisting of some ``nodes'' and some ``links'' 
(or edges). When
each pair of node is connected by a link, the network becomes a trivial
one. One therefore links (or does not link) two nodes 
according to some intrinsic property of them. Depending on the context, 
different properties of the nodes are relevant in different networks, leading
to different rules for linking. For example, in science collaboration network,
each scientist is a node and two nodes (scientists) are linked when they
are co-authors in at least one paper. In English language network, each word
is a node and a pair of words are linked if (in one or more sentences) they
appear side by side or one word apart. In this way, a network structure can 
be identified in widely varied contexts. 
Once a network is identified, one can measure in it some characteristic 
properties like the average shortest distance, clustering coefficient, degree 
distribution etc. When the degree distribution 
decays as a power law, the network is said
to be scale-free. When the average shortest distance is small (increases only
logarithmically with the size of the network) but the clustering coefficient 
is high compared to the random network, the network is 
called small world. Many natural and man-made networks \cite{review} have 
been found to be scale-free and/or small-world.

Recently, Corso \cite{Corso} has considered a network where each node is a
natural number and two nodes are linked if they share a common prime factor
larger than some chosen lower limit. The network is not scale-free (except in
some restricted sense) but is a 
small-world unless the lower limit is 1, that is, unless one considers all 
prime factors.
Motivated by this study, we consider here a network where each node is a prime
number. The rule of placing links will be explained in the next section. The
rule relies upon the validity of what is known as Goldbach conjecture 
\cite{Schroder} and involves a tunable parameter. Depending on the value 
of the parameter, we have a small-world or a regular network. 
We mention that our work has no connection with the 
issue of the {\em validity} of Goldbach conjecture.

In the next section, we shall describe the network and present a theoretical 
analysis of its behaviour. In Section III  we shall describe the 
computational studies and in Section IV present the general conclusions.

\section{The Model}

Goldbach conjecture says that any even number ($>2$) can be written as the sum 
of two prime numbers (often in more than one way). To construct our network,
we start with the even number 8 and note that it can be broken up into primes
as $8=3+5$. (For avoiding uninteresting complications we do not consider even
numbers below 8.) We put (the first) two nodes in the network and label them 
as 3 and 5 and put a link between them. Now we consider all even numbers 
10, 12, 14, $\cdots$ upto some $N_e$ and break up each of them into primes as
$n = p + q$. If $p$ is not already in the network, a new node labelled as $p$ 
is added, and similarly for $q$. Then a link is put between $p$ and $q$. A
complication is that (as mentioned earlier) very often one even number can be
broken up into primes in more than one way. If one puts links corresponding
to all the decompositions, then one has a link between almost every pair of 
nodes and the network becomes trivial. Therefore, with 
the following prescription we choose {\em one} way of decomposition 
for every even number, depending on the {\em difference} between the component 
primes.  We calculate the difference $\Delta =
|p - q|$ between the two components for every break-up and choose one 
break-up with probability ${\Delta}^{\alpha}$, where $\alpha$ is a (in
fact, the only) parameter of the model. 
For example, the number $n=24$ can be broken up in three ways : 5+19, 7+17, 
and 11+13, with $\Delta$=14, 10, and 2 respectively. In a large number of 
realisations, the link for 24 will be between 5 and 19 with probability 
$p_1$, between 7 and 17 with probability $p_2$, and between 11 and 13 
with probability $p_3$, where $p_1 = 14^{\alpha}/s$, $p_2 = 10^{\alpha}/s$
and $p_3 = 2^{\alpha}/s$ with $s = 14^{\alpha} + 10^{\alpha} + 2^{\alpha}$.
As $p_1 + p_2 + p_3 = 1$, we could realise the choice of prime-pair 
by calling a random number between 0 and 1. 

One should note that since we put one link for each even number, 
for $M$ even numbers one will have exactly $M$ links,
but some $N$ ($<M$) number of nodes. One should also note that through
the parameter $\alpha$ we actually control the difference between the prime
pairs (chosen to be linked) for each even number. Thus, for $\alpha=0$, 
our choice is independent of the difference $\Delta$,
while for $\alpha = - \infty $ ($+ \infty$) the break-up with the 
smallest (largest) difference between the components is chosen.

What type of network we have thus constructed ? To have an answer
analytically, let us calculate the average value of $\Delta$ for a given
even number $n$. This is
\[
<\Delta> = \frac{1}{\Omega} \sum {\Delta}^{\alpha + 1},
\]
where the sum extends over all Goldbach pairs of $n$ and $\Omega$ is the
number of such pairs. For positive
values of the power $(\alpha + 1)$, the value of this sum will be dominated
by the terms with large $\Delta$ and $<\Delta>$ will hence increase as $n$
increases. One of the chosen Goldbach pair will then become small, and these
small numbers will be highly populated nodes. The network may therefore be
expected to be of small-world type.
On the other hand, for negative values of $(\alpha + 1)$, the value of the
sum will be dominated by small $\Delta$ values and for large $n$ the quantity
$<\Delta>$ will converge to some finite number. Both the chosen Goldbach
pairs will be large (so that the difference between them remains small)
and no node will be very highly populated. The network is then likely to
lose the small-world character. Although $\Delta$ runs over some (but not
all) of the natural numbers, the convergence behaviour of $<\Delta>$ may be
expected to be the same as the Riemann zeta function $\zeta (-\alpha - 1)$.
As this function converges only for $\alpha < -2$, the change-over in the
behaviour of the network is expected to occur around $\alpha = -2$.

\section{Simulation Studies }

To analyse the properties of the network by computer simulation, we measure
several characteristics of the network.
(i) {\em Average shortest distance between two nodes (d) } : The 
shortest distance between two nodes is the smallest number of links via which
one can go from one node to the other. We have measured this quantity for
all pairs of nodes and taken the average. Results for the 
measurement of this quantity is presented in Fig. 1 as a function of the
number of nodes ($N$). It is observed that this quantity varies linearly
with logarithm of the number of nodes as it happens for a small-world
network. This behaviour prevails for all values 
of $\alpha$ upto a lower limit of $\alpha_0 = -1.8$.  In particular, 
as $\alpha$ varies from 5 to 1, the $d$ - $N$ (log-linear) plot moves upwards 
parallel to itself. As $\alpha$ decreases further, the lines
continue to move upwards but the slope increases continuously until for
$\alpha < \alpha_0$ the line ceases to be straight and starts bending upwards.
In this region, $d$ varies linearly with $N$, as it happens for a regular 
network. 

\begin{figure}
\noindent \includegraphics[width= 8cm, angle = 270]{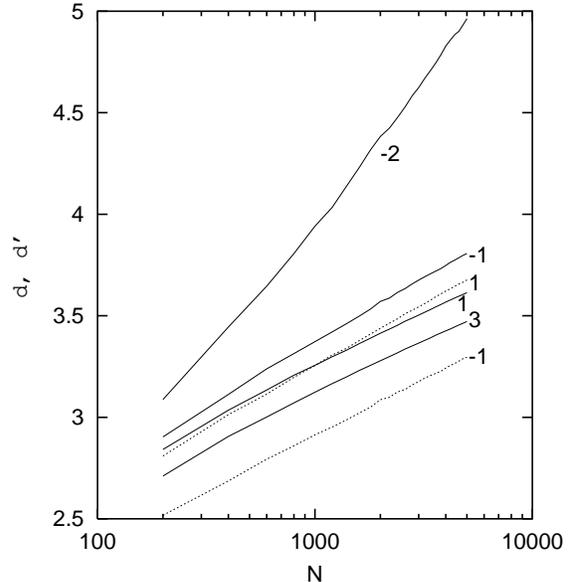}
\caption{Average shortest distance ($d$) as a function of the
number of nodes
for the network constructed from prime numbers (continuous line). Also
shown is the same plot for a random network having the same number
of nodes ($d^{\prime}$, dotted line).
The lines for $d$ and $d^{\prime}$ cross over at $\alpha=1$ and $N=1000$.
The numbers labelling the curves stand for the value of $\alpha$. All results
presented in this paper have been averaged over about 20 realisations of the
network.}
 \end{figure}

\begin{figure}
\noindent \includegraphics[width= 8cm, angle = 270]{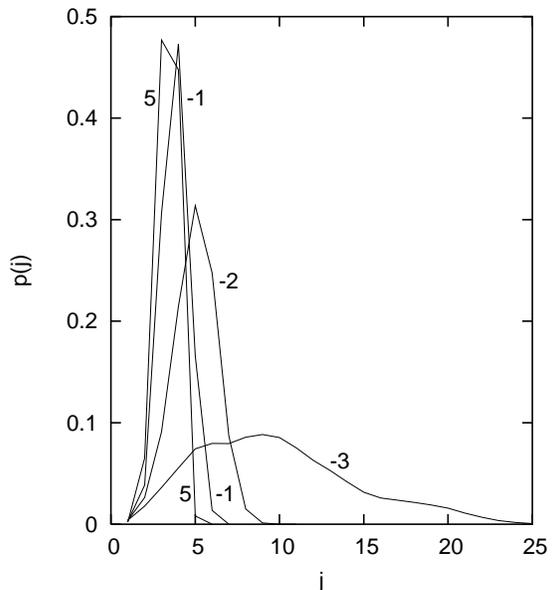}
\caption{ Plot of $p(j)$ as a function of $j$, where $p(j)$ is the
probability that a pair of nodes chosen randomly from the network will be
$j$ distance apart. The numbers labelling the curves stand for the
value of $\alpha$. $N=5000$.}
 \end{figure}

We have also plotted in Fig. 1, the average shortest distance
($d^{\prime}$) for a random network having the same number of nodes. As is
well-known \cite{review}, the $d^{\prime}$ - $N$ (log-linear) plot is a
straight line for the entire range of values of $\alpha$ and the lines
maintain a constant slope and move gradually upwards as $\alpha$ increases.
(However, for $\alpha > 2$ the $d^{\prime}$ - $N$ plot does not depend much
on the value of $\alpha$.)

\begin{figure}
\noindent \includegraphics[clip,width= 8cm, angle = 270]{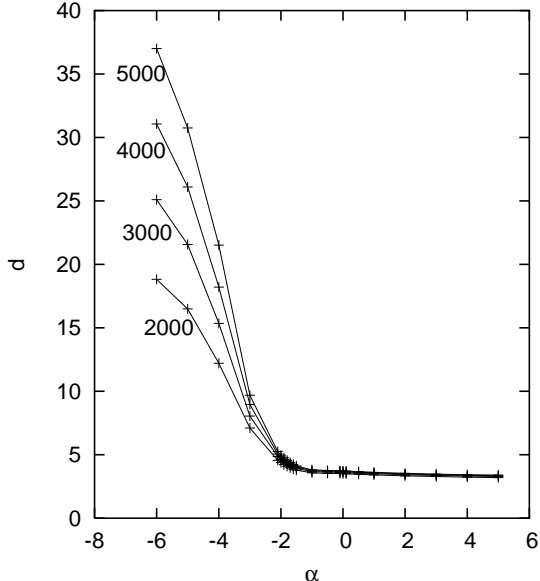}
\caption{ Plot of average shortest distance $d$ as a function of $\alpha$
for different values of the number of nodes $N$. The numbers labelling the
curves stand for the value of $N$.}
 \end{figure}

To gain further insight into the change of behaviour of $d$ as a function of
$\alpha$, we have measured (Fig. 2) the probability $p(j)$ that a pair of 
nodes 
chosen randomly from the network will be $j$ distance apart. It is observed 
that for $\alpha > \alpha_0$, $p(j)$ is high only for small values of $j$, 
indicating that most of the pairs of nodes are at a small distance apart
and $d$ (which is nothing but $\sum_j j p(j)$) is small. On the other hand,
for $\alpha < \alpha_0$, the distribution $p(j)$ is flat over a large range 
of values of $j$, 
indicating that the distance between a pair of nodes will also often be 
large and $d$ will be high due to the contribution from high $j$-values. 

Lastly, in order to ascertain how sharply the change of behaviour occurs at
$\alpha=\alpha_0$, we plot $d$ against $\alpha$ for different values of
$N$ (Fig. 3). The rise of $d$ for $\alpha < \alpha_0$ becomes sharper and
sharper as $N$ increases. Our estimate of
$\alpha_0 = -1.8$ is based on simulations of networks with at most 5000 nodes.

(ii) {\em Clustering coefficient} ($C$) : The clustering coefficient $C_i$ for 
the node $i$, is defined as the ratio $M_i/m_i$ where $M_i$ is the actual 
number of links among the {\em neighbours} of the node $i$, and $m_i$ is the 
number of all possible links among the neighbours of $i$. (Thus, if $i$ has
degree $k_i$, then $m_i=k_i(k_i+1)/2$.) The clustering coefficient for the
entire network is defined as the average of $C_i$ over all nodes $i$.
This quantity $C$ has been measured in
our network and compared with the same ($C^{\prime}$) for a random network 
with the same number of nodes (Fig. 4). It is found that $ C > C^{\prime}$ 
for the entire range of values of $\alpha$ investigated here. This behaviour,
combined with the results for $d$, leads us to conclude that the network is of
{\em small-world type for} $\alpha > \alpha_0$ and of {\em regular type
for} $\alpha < \alpha_0$.

\begin{figure}
\noindent \includegraphics[width= 8cm, angle = 270]{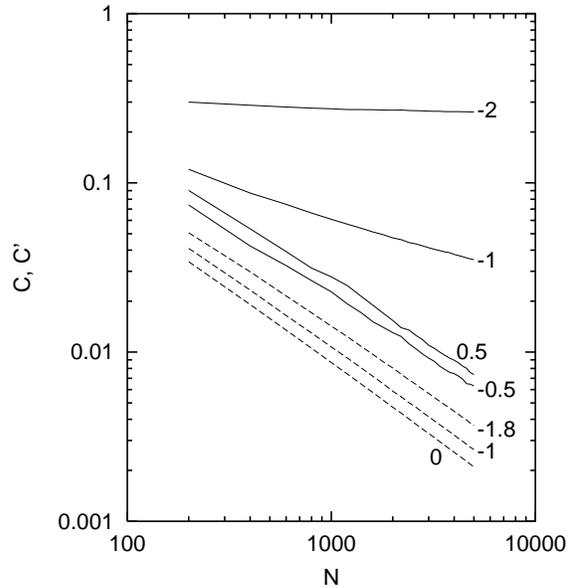}
\caption{Clustering coefficient ($C$) as a function of the number
of nodes ($N$) for the network constructed from prime numbers (continuous
line)i. Also shown is the clustering coefficient $C^{\prime}$ (dotted line)
for for a random network having the same number of nodes.
For the entire range of $\alpha$ values $C$ remains larger than $C^{\prime}$.
The numbers labelling the curves stand for the value of $\alpha$.}
\end{figure}

The details of the behaviour of the clustering coefficient as a function of 
the number of nodes is as follows. For a given $\alpha$, $C$ decays 
algebraically with $N$ and the $C$~-~$N$ line rises, maintaining a constant 
slope as $\alpha$ is increased from -0.5.
But as $\alpha$ is decreased from -0.5, the line remains straight,
rises upwards but becomes more and more horizontal. The slope almost
vanishes (particularly in log-log scale) for $\alpha < \alpha_0$. For the 
corresponding random network, the
clustering coefficient $C^{\prime}$ also decreases algebraically with $N$, and
the $C^{\prime}$ - $N$ line rises,
maintaining a constant slope, as $\alpha$ is decreased from 0 to $\alpha_0$.
The lines for $\alpha < \alpha_0$ are almost
coincident with those of $\alpha = \alpha_0$ and the lines for $\alpha > 0$
are also almost coincident with those of $\alpha = 0$.

(iii) {\em The degree distribution function} $P(k)$ (defined as 
the probability that a node has $k$ links attached to it) is of irregular 
type (Fig. 5) 
for all reasonable values of $\alpha$, indicating that the network is not of 
scale-free nature. However, some other related plots do contain some
interesting information. Thus, it is interesting to observe how the network 
grows for different values of $\alpha$ (Fig. 6). For large positive values 
of $\alpha$, by breaking up an even number $n$ one chooses only those pairs
of primes that are far apart. One of the chosen primes will therefore be a
small prime number, while the other will be one that is close to $n$. Very
often one will find that the latter prime has not been included in the
network till now and thus a new node is added. The network then grows in size.
On the other hand, for small (large negative) values of $\alpha$, the
difference between the primes chosen will be small. Each prime will then be
$\sim (n/2)$ and will very often be found to be already present in the network.
The network will then grow slowly. Such behaviour has been confirmed by
simulation (Fig. 6).

\begin{figure}
\noindent \includegraphics[width= 8cm, angle = 270]{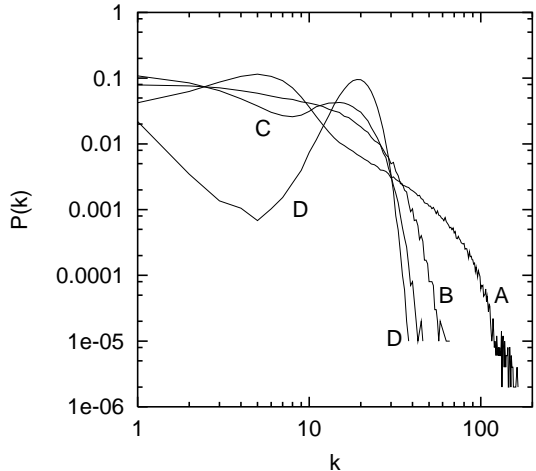}
\caption{ Degree distribution $P(k)$ in the network constructed from
prime numbers. The different lines correspond to $\alpha =$ 2 (A), -0.1 (B),
-0.5 (C), -2 (D). $N=5000$.}
 \end{figure}

We have also studied the $\alpha$ dependence (Fig. 7) of average connectivity 
$< k >$ (which is simply $\sum_k k P(k) = 2M/N$,  $M$ being number of links  and
$N$ being number of nodes) and fluctuation in connectivity defined as
\[
f(k) = \surd (<k^2> - <k>^2).
\]
Both these quantities display a sharp change at $\alpha = \alpha_0$. For
a given number of nodes, the value of $<k>$ will be large for 
$\alpha < \alpha_0$ and small for $\alpha >  \alpha_0$ since there are more 
links in the former case than in the latter.
For a regular network, the nodes have
the same degree rather uniformly, so that $f(k)$ 
is small, but for a small-world network, some nodes
are very rich in degree while the other nodes have low degree, rendering
$f(k)$ very large. Moreover, for small-world network ($\alpha > \alpha_0$)
the rich nodes go on gaining links as the network evolves, so that
the degree of the most-connected node ($k_m$, say) increases with size
of the network (Fig. 8). In the regular network regime, no node is 
preferentially linked and $k_m$ does not increase very much with $N$ and
in fact approaches the average connectivity $< k >$.

\begin{figure}
\noindent \includegraphics[width= 8cm, angle = 270]{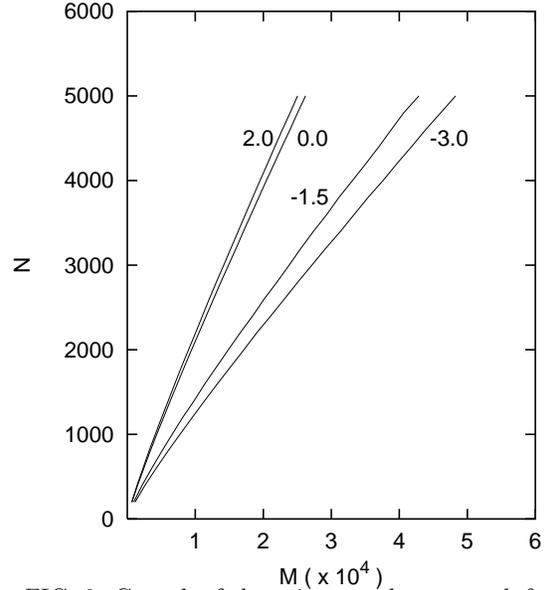}
\caption{ Growth of the prime number network for different values of 
$\alpha$. The number of nodes ($N$) is plotted as a function
of the number of links ($M$). Since at each time step one link is added, the
X-axis also represents the time step. The numbers labelling the curves stand 
for the value of $\alpha$.}
 \end{figure}
 
\begin{figure}
\noindent \includegraphics[width= 8cm, angle = 270]{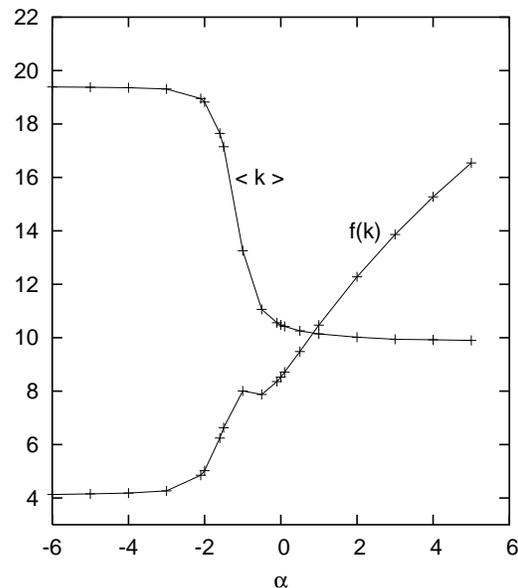}
\caption{ Average connectivity $<k>$ and fluctuation in connectivity $f(k)$
as a function of $\alpha$. $N=5000$. }
 \end{figure}

\begin{figure}
\noindent \includegraphics[width= 8cm, angle = 270]{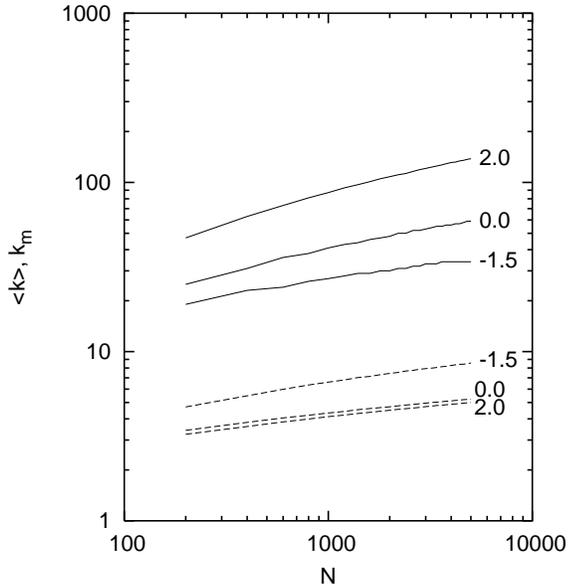}
\caption{ Average connectivity $<k>$ (broken line) and connectivity of the best 
connected node $k_m$ (continuous line) as a function of the number of 
nodes ($N$) for different
values of $\alpha$. As $\alpha$ decreases, the $<k>~-~N$ line rises, while the
$k_m~-~N$ line descends. The numbers labelling the curves stand
for the value of $\alpha$.}
 \end{figure}

(iv) {\em Clustering coefficient $C(k)$ measured as a function of the degree} 
of a node has been proved to be an important characteristic for many real 
networks \cite{review,rail}. This quantity is defined as $C_i$ (as defined
above) averaged by running $i$ over only those nodes that have degree $k$. 
For $\alpha < 0$,
this quantity has a peak at a small ($< 20$) value of $k$, indicating that
the nodes mostly have low degree and high clustering coefficient. On the other
hand, for $\alpha > 0$, $C(k)$ is almost flat extending over a large range of
values of $k$. (Fig. 9)

(v) {\em Degree-degree correlation function} $r$, as proposed by Newman
\cite{Newman} measures the tendency of a link to have same type of degree
(both high or both low) at its two ends. Thus, when $r$ is positive, the
network is assortative, and a link prefers to have the same type of node at
the two ends whereas, when
$r$ is negative, the network is disassortative, and a link prefers to have
different type of nodes (one of high degree and the other of low degree) at the 
two ends. This parameter may be measured from the relation \cite{Newman}
\[
r = \frac{M^{-1}\sum_i j_i k_i - [M^{-1}\sum_i (j_i + k_i)/2]^2}
{M^{-1}\sum_i (j_i^2 + k_i^2)/2 - [M^{-1}\sum_i (j_i + k_i)/2]^2}
\]
where $j_i$ and $k_i$ are the degrees of the nodes that are at the two ends
of the $i$-th link.
For the network under study, the quantity $r$ has been found to be
positive (negative)
for negative (positive) values of $\alpha$ (Fig. 10). This indicates that the 
nature of the network changes from assortative to disassortative as 
$\alpha$ changes its sign. 

\begin{figure}
\noindent \includegraphics[clip,width= 8cm, angle = 270]{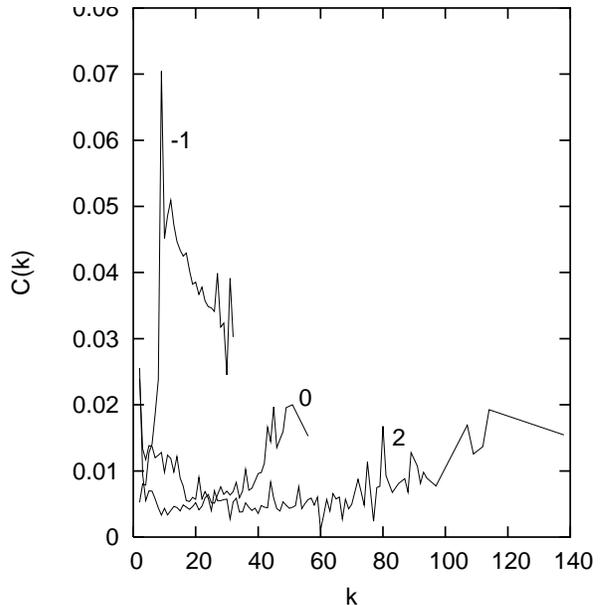}
\caption{Clustering coefficient $C(k)$ measured as a function of
the degree. The numbers labelling the curves stand for the 
value of $\alpha$. For $\alpha < -1$ the curve looks similar to that
for $\alpha = -1$ but goes higher up and for $\alpha > 2$ the curve becomes 
more flat and moves further down. $N=5000$.}
 \end{figure}

\begin{figure}
\noindent \includegraphics[clip,width= 8cm, angle = 270]{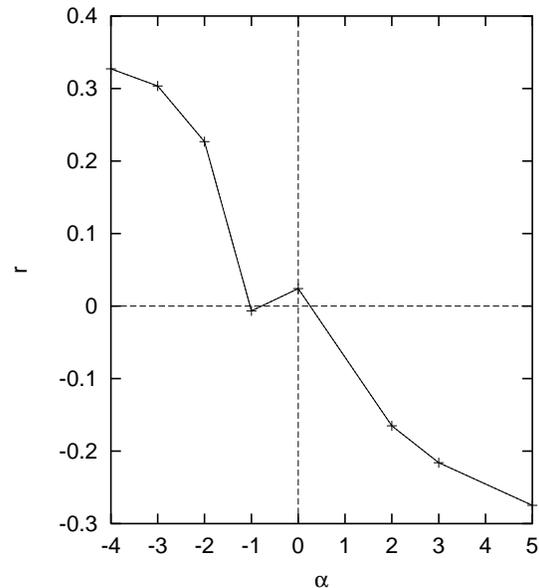}
\caption{Degree-degree correlation coefficient $r$ as a function of
$\alpha$. Note that $r$ bears a sign opposite to that of $\alpha$. $N=5000$.}
 \end{figure}

\section{Conclusion}

In conclusion, we have constructed a network of prime numbers with links
placed on the basis of Goldbach conjecture. The network is of small world type
when a parameter $\alpha$ of the model is larger than $-1.8$ and of regular 
type when $\alpha$ is lower than $-1.8$. One must note that,
for $\alpha > 0$ larger values of $\Delta $ (difference between the 
component primes) are preferred and the addition of a
new link leads to a large prime getting attached to a small one. Thus,
preferrential attachment to small primes are supported and the `rich gets
richer' principle leads to a small-world network, as for the case
of Barabasi-Albert network [1] although the scale-free property is not
observed. In any case, the small-world property indicates 
some pattern in the Goldbach decomposition of prime numbers.

{\bf Acknowledgement :} The authors are grateful to Marc Barthelemy and
Sitabhra Sinha for encouragement and fruitful discussions. The work of one 
author (AKC) was supported by UGC fellowship.

\end{multicols}

\end{document}